\documentclass[12pt]{article}

\usepackage{amssymb}
\usepackage{amsmath}

\newcommand{\Lg}{\mbox{$\mathfrak g$}}
\newcommand{\Ll}{\mbox{$\mathfrak l$}}
\newcommand{\Lh}{\mbox{$\mathfrak h$}}
\newcommand{\Lk}{\mbox{$\mathfrak k$}}
\newcommand{\Lp}{\mbox{$\mathfrak p$}}

\newcommand{\Lm}{\mbox{$\mathfrak m$}}

\newcommand{\Lq}{\mbox{$\mathfrak q$}}

\newcommand{\Pf}{{\em Proof}. }
\newcommand{\EPf}{\hfill$\square$}

\newcommand{\Ad}[1]{\mbox{$\mbox{Ad}_{#1}$}}

\newtheorem{thm}{Theorem}

\newtheorem{lem}[thm]{Lemma}

\begin{document}

\title{Variationally complete actions\\ on compact symmetric spaces}
\author{Claudio Gorodski\footnote{Partially supported by CNPq and
FAPESP.}\hskip .3cm and Gudlaugur Thorbergsson}

\maketitle

\begin{abstract}
We prove that an isometric action of a compact Lie group on a
compact symmetric space is variationally complete if
and only if it is hyperpolar.
\end{abstract}

\section{Introduction}

The main result of this paper is the following theorem.

\medskip\noindent\textbf{Theorem}
\textit{An isometric action of a compact Lie group on a compact symmetric
space is variationally complete if
and only if it is hyperpolar.}\medskip

Variationally complete actions were
introduced by Bott in \cite{B} (see also~\cite{BS}).
Let $G$ be a compact Lie group acting on a complete
Riemannian manifold $M$ by
isometries. A geodesic $\gamma$ in $M$  is called
\emph{$G$-transversal} if it is orthogonal to the $G$-orbit through
$\gamma(t)$
for every $t$. One can show that a geodesic
$\gamma$ is $G$-transversal if there is a point $t_0$ such that
$\dot\gamma(t_0)$ is orthogonal to $G\gamma(t_0)$.
A Jacobi field along a
geodesic in $M$ is called \emph{$G$-transversal}
if it is the variational
vector field of a
variation through
$G$-transversal geodesics. The action of
$G$ on $M$ is called
\emph{variationally complete} if every $G$-transversal
Jacobi
field $J$ in $M$ that is tangent to the $G$-orbits at two
different
parameter values is the restriction of a Killing
field on
$M$
induced by the $G$-action.
It is proved in~\cite{BS} on
p.~974 that instead
of requiring tangency at two
different points in the definition of
variational completeness it is
equivalent to require tangency at one point
and
vanishing at another point.

Conlon considered in~\cite{C} actions of a
Lie group $G$ on a complete
Riemannian manifold $M$ with the property that
there is a
connected submanifold $\Sigma$ of $M$ that meets
all orbits of
$G$ in such a way that the intersections between $\Sigma$ and
the orbits of
$G$ are all orthogonal. Such a submanifold
is called  a \emph{section} and
an action admitting a section is
called \emph{polar}.
Notice that neither do
we assume as Conlon in~\cite{C} that $\Sigma$ is
closed nor do
we assume
that it is properly embedded as is usually required in the recent
literature
on the subject.
It is easy to see that a section
$\Sigma$ is totally
geodesic in
$M$. An action admitting a section that is flat in the induced
metric is
called \emph{hyperpolar}.
Conlon proved in~\cite{C} that
a
hyperpolar action of a compact Lie
group on a complete
Riemannian manifold
is variationally complete. Notice that he does not use
in his proof that
the
sections of the action are closed. His result therefore implies
one
direction of the main theorem
of this paper.

In the case of Euclidean
spaces,
hyperpolar representations were classified by Dadok in~\cite{D}.
As
a consequence of his classification he obtained that
a hyperpolar
representation of a compact
Lie group is orbit equivalent to the
isotropy
representation of a symmetric space.
(We recall
that two isometric
actions are said to be
\emph{orbit equivalent} if there is an
isometry
between the action spaces under which the orbits of
the two actions
correspond.)

On the other hand, we classified variationally
complete
representations in~\cite{GTh}. As a consequence
of our
classification we obtained that a variationally complete
representation of a
compact Lie group is hyperpolar.
Di Scala and Olmos gave in~\cite{DO} a very
short, simple proof
of this result. It follows that a representation is
hyperpolar
if and only if it is variationally complete.

Next we discuss the
case of compact symmetric spaces.
There are two important classes of
examples of hyperpolar
actions on them. It is clear that a cohomogeneity
one
action of a compact Lie group on a compact symmetric space is
hyperpolar.
Hermann constructed
in~\cite{H} a class
of examples of variationally
complete actions on compact symmetric spaces
which in fact turned out to be
hyperpolar. Namely, if $K_1$ and $K_2$
are two symmetric subgroups of the
same compact Lie group $G$, then the
action of $K_1$ on $G/K_2$ is
hyperpolar and so is the action of
$K_1\times K_2$ on $G$. Kollross~\cite{K}
classified hyperpolar actions
on compact \emph{irreducible}
symmetric
spaces. It follows from his classification that,
in the irreducible case,
all hyperpolar actions belong to either one of
these
classes.

Finally, we
add that the concepts of hyperpolarity
and variational completeness admit
natural extensions
in the context of proper Fredholm actions of Hilbert-Lie
groups
on Hilbert spaces (see~\cite{T1} and \cite{T2}).
We shall make use of
them in the course of our proof.

\section{The proof of the theorem}

Conlon
proved in~\cite{C} that hyperpolar
actions are variationally complete as we
pointed out in the introduction.
In this section we will prove the converse
of his result for actions
on compact symmetric spaces. Our strategy is as
follows.
First we reduce to the case of a symmetric space
of compact type.
Next we lift the action on the symmetric space to a
variationally complete
action of a path group on a Hilbert space.
Then we use an argument similar
to that in~\cite{DO}
to show that the lifted action is hyperpolar. Finally,
it is easy to see
that a section for the action on the Hilbert space induces
a
flat section for the original action.

Let $M$ be a compact Riemannian
symmetric space. We identify $M$
with the coset space $G/K$, where $G$ is
the connected
component of the group of all isometries of $M$ and $K$ is the
isotropy
subgroup  of a chosen base point. Let $H$ be a compact Lie group
acting by isometries on $M$.
An action of a group
is variationally complete (resp.~hyperpolar) if and only if
the same is true for the restriction of that action to the connected
component of the group. This follows from the definition 
in the case of variationally complete actions and is easy in the
case of hyperpolar actions, see Proposition~2.4 in~\cite{HPTT2}.
We can thus assume that $H$ is a connected,
closed subgroup of $G$.

Let $\hat M$ be a compact Riemannian covering space of $M$ which splits as
a product of a torus $T^k$ and a symmetric space $N$ of compact type.
Let $\hat G$ be the
connected component of the group of all isometries of $\hat M$. Then $\hat
G$ is a covering group of $G$ and we let $\hat H$ denote the
subgroup of $\hat G$ which is the connected
component of the inverse image of $H$.
Notice that the action of $\hat H$ on $\hat M$ is
variationally complete if
the one of $H$ on $M$ is so. In this case, the restriction
$\sigma_2$ of the action of $\hat H$ to $N$ is
clearly still variationally complete.
We will prove below that a variationally complete
action on a symmetric space of compact type is hyperpolar.
Since the restriction $\sigma_1$ of the action of
$\hat H$ to $T^k$ is orbit equivalent to the action of a torus on $T^k$ and
therefore hyperpolar, it follows that the product action
$\sigma_1\times\sigma_2$ on $T^k\times N=\hat M$ is
hyperpolar. The action
of $\hat H$ on $M$ is therefore also hyperpolar since it is
orbit equivalent
to the product action
of $\sigma_1\times\sigma_2$ on $\hat M$,
which can be seen by arguments
similar to those in the proof of Theorem~4~(ii) in~\cite{D} or
of Proposition~3.4~(d) in~\cite{GTh}.
Now we can simply project a section of the
$\hat H$-action in 
$\hat M$ down to $M$ to see that the action of $H$ on $M$ is hyperpolar.

It follows from the discussion in the previous paragraph
that from now on we can restrict to the case where $M$ is of compact type.
In this case the group $G$ is semisimple.
There is an involution $\sigma$ of $G$
such that $K$ is open in the fixed point set of $\sigma$.
Now the Lie algebra of $G$ decomposes into the $\pm1$-eigenspaces
of $d\sigma$, $\Lg=\Lk+\Lp$, where $\Lk$ is the Lie algebra
of $K$, and the Riemannian metric in $M$ is the $G$-invariant
metric induced by some $\Ad{G}$-invariant inner
product on $\Lg$. The group $H$ is a closed subgroup of $G$
so that it acts on $M$ by left
translations, and $H\times K$ acts on $G$ by $(h,k)\cdot g= hgk^{-1}$,
where $h\in H$, $k\in K$ and $g\in G$. Clearly, the projection
$\pi:G\to G/K$ is an equivariant Riemannian submersion.

\begin{lem}
If the action of $H$ on $M$ is variationally complete,
then the action of $H\times K$ on $G$ is also variationally complete.
\end{lem}

\Pf Let $\gamma$ be an $H\times K$-transversal geodesic in~$G$ defined
on~$[0,1]$ and let $J$ be an $H\times K$-transversal Jacobi field
along $\gamma$ such that $J(0)=0$ and $J(1)$ is tangent to the
$H\times K$-orbit through $\gamma(1)$. We must show that
$J$ is the restriction along $\gamma$ of an $H\times K$-Killing field.

We
first identify the tangent and normal spaces to the
$H\times K$-orbits in
$G$. Let $\Lh$ be the Lie algebra of $H$.
Then $\Lh$ is a subalgebra of
$\Lg$ and we denote by
$\Lh^\perp$ its orthogonal complement in~$\Lg$ with
respect to the
$\Ad G$-invariant inner product on $\Lg$. Note
that
$[\Lh,\Lh^\perp]\subset\Lh^\perp$
and that $\Lk$ and $\Lp$
are also
mutually orthogonal in $\Lg$.
Given $a\in G$, the tangent space
to $H\cdot
a\cdot K$ at the point $a$ is $\Lh\cdot a + a\cdot\Lk$.
Therefore the left
translates under $a^{-1}$ of the tangent and normal
spaces of $H\cdot a\cdot
K$ at $a$ are respectively
$\Ll^a:=\Lk+\Ad{a^{-1}}\Lh$ and
$\Lm^a:=\Lp\cap\Ad{a^{-1}}\Lh^\perp$.
It will be convenient to set
$\Lq^a:=\Lk\cap\Ad{a^{-1}}\Lh$.

Let $\gamma(0)=a\in G$. Then
$\gamma(t)=ae^{tX}$ for
$X\in\Lm^a$. Since $J(0)=0$, setting
$\gamma_\alpha(t)=ae^{t(X+\alpha Y)}$,
where $Y\in\Lm^a$ and $Y=a^{-1}\cdot
J'(0)$,
defines a variation of $\gamma=\gamma_0$ through
$H\times
K$-transversal geodesics which induces $J$.
Since $\pi:G\to G/K$ is an
equivariant Riemannian submersion and
each $\gamma_\alpha$ is a horizontal
curve with respect to $\pi$,
we have that $\{\pi\gamma_\alpha\}$ is a
variation of
$\bar\gamma=\pi\gamma$ through
$H$-transversal geodesics in
$M$. Moreover, the associated $H$-transversal
Jacobi field $\bar J$ along
$\bar\gamma$ satisfies $\bar J(0)=0$,
$\bar J'(0)=a\cdot\pi_*Y$,
and $\bar
J(1)$ is tangent to the $H$-orbit
through $\pi(ae^X)$.
It follows by
variational completeness of the
action of $H$ on $M$ that $\bar J$ is the
restriction along $\bar\gamma$ of
an $H$-Killing field on $M$, namely $\bar
J(t) = Z\cdot \bar\gamma(t)$
for some $Z\in\Lh$, which yields
\[
 \bar J(t)
= \frac{\partial}{\partial\alpha}\Big\arrowvert_{\alpha=0}
e^{\alpha
Z}\pi(ae^{tX}) 
=
\frac{\partial}{\partial\alpha}\Big\arrowvert_{\alpha=0}ae^{tX}\pi(e^{-tX}
a^
{-1}e^{\alpha Z}ae^{tX})
=ae^{tX}\cdot\pi_*\Ad{(ae^{tX})^{-1}}Z.
\]
This
equation gives
$\bar J(0)=a\cdot\pi_*\Ad{a^{-1}}Z$. Since we already know
that
$\bar J(0)=0$, we get that $\pi_*\Ad{a^{-1}}Z=0$, and
then
$\Ad{a^{-1}}Z\in\Lk$. Hence $\Ad{a^{-1}}Z\in\Lq^a$.

Since $X\in\Lp$,
the left translations by $e^{tX}$ define isometries
of the symmetric
space~$M$ which form a one-parameter group of
transvections, that is,
each
one of them induces parallel transport along the
geodesic
$s\mapsto\pi(e^{sX})$. Denote by $\nabla$ the
Levi-Civita
connection of $M$. Then we have
\[ a^{-1}\cdot\bar
J'(0)=a^{-1}\cdot\nabla_{a\cdot\pi_*X}\bar
J
=\nabla_{\pi_*X}(a^{-1}\cdot\bar
J)=\frac{d}{dt}\Big|_{t=0}\pi_*\Ad{(ae^{t
X})^{-1}}Z=\pi_*[\Ad{a^{-1}}Z,X].
\]

We deduce that $\bar
J'(0)=a\cdot\pi_*[\Ad{a^{-1}}Z,X]$.
But we already know that $\bar
J'(0)=a\cdot\pi_*Y$.
Since $Y\in\Lm^a\subset\Lp$
and
$[\Ad{a^{-1}}Z,X]\in[\Lq^a,\Lm^a]\subset\Lm^a\subset\Lp$, this
implies
$Y=[\Ad{a^{-1}}Z,X]$.

Finally, consider the $H\times K$-Killing
field in $G$ given
by 
\[ U\cdot y =
\frac{\partial}{\partial\alpha}\Big\arrowvert_{\alpha=0}
e^{\alpha
Z}ye^{-\alpha\Ad{a^{-1}}Z}. \]
Then
\[
U\cdot\gamma(t)=\frac{\partial}{\partial\alpha}\Big\arrowvert_{\alpha=0}
e^{
\alpha Z}ae^{tX}a^{-1}e^{-\alpha Z}a. \]
Therefore $U\cdot\gamma(0)=0$
and
\begin{eqnarray*}
\frac{d}{dt}\Big\arrowvert_{t=0}U\cdot\gamma(t)&=&\frac{\partial}
{\partial\alpha}\Big\arrowvert_{\alpha=0}
\frac{\partial}{\partial
t}\Big\arrowvert_{t=0}e^{\alpha
Z}ae^{tX}a^{-1}e^{-\alpha Z}a\\
&=&(\frac{\partial}{\partial\alpha}\Big\arrowvert_{\alpha=0}
\Ad{e^{\alpha
Z}a}X)\cdot a\\
 &=& [Z,\Ad{a}X]\cdot a.
\end{eqnarray*}
Hence
$a^{-1}\cdot\frac{d}{dt}\big\arrowvert_{t=0}U\cdot\gamma(t)=
\Ad{a^{-1}}[Z,\Ad{a}X]=Y.$
It follows that $J(t)=U\cdot\gamma(t)$. \EPf

\medskip

Let
$V=L^2([0,1];\Lg)$ denote the Hilbert space of $L^2$-integrable
paths
$u:[0,1]\to\Lg$, and let $\hat G=H^1([0,1];G)$ denote the
Hilbert-Lie
group of $H^1$-paths in $G$ parametrized on~$[0,1]$ (the
elements
of $\hat
G$ are the absolutely continuous paths $g:[0,1]\to G$
whose derivative is
square integrable).
We have that $\hat G$ acts by affine isometries
on $V$:
$g*u=gug^{-1}-g'g^{-1}$, where $g\in\hat G$ and $u\in V$.
Let $\mathcal
P(G,H\times K)$ denote the subgroup of all paths
$g\in\hat G$ such that
$(g(0),g(1))\in H\times K$.
Let $\varphi:V\to G$ be the
parallel transport
map defined by $\varphi(u)=h(1)$, where
$h\in\hat G$ is the unique solution
of $h^{-1}h'=u$, $h(0)=1$.
Then it is known that (see~\cite{T1,T2,
TTh}):
\begin{enumerate}
\item[(a)] the action of $\hat G$ on $V$ is proper
and Fredholm;
\item[(b)] $\varphi(g*u)=g(0)\varphi(u)g(1)^{-1}$, where
$g\in\hat G$ and
$u\in V$;
\item[(c)] $\mathcal P(G,H\times
K)*u=
\varphi^{-1}((H\times K)\cdot\varphi(u))$;
\item[(d)] the action of
$\mathcal P(G,1\times G)$ on $V$ is simply
transitive;
\item[(e)]
$\varphi:V\to G$ is a Riemannian submersion and a
principal
$\Omega_1(G):=\mathcal P(G,1\times 1)$-bundle;
\item[(f)] the horizontal
distribution $\mathcal H$ of $\varphi$ is given
by $\mathcal
H(u)=\{\Ad{g}\hat Y:Y\in \Lg\}$, where $\hat Y$ denotes the
constant
path
with value $Y$ and $g$ is the unique element
of $\mathcal P(G,1\times
G)$ which satisfies $u=g*\hat 0$;
\item[(g)] the tangent space to the orbit
through $\hat Y$
at $\hat Y$ is 
\[ \{[\xi,\hat Y]-\xi':\;\xi\in
H^1([0,1],\Lg),\;
\xi(0)\in\Lh,\; \xi(1)\in\Lk\}.\]
\end{enumerate}

\begin{lem}
If the action of $H\times K$ on $G$ is variationally
complete,
then the action of $\mathcal P(G,H\times K)$ on $V$ is
also
variationally complete.
\end{lem}

\Pf Let $\gamma$ be a transversal
geodesic in~$V$ defined
on~$[0,1]$ and let $J$ be a transversal Jacobi
field
along $\gamma$ such that $J(0)=0$ and $J(1)$ is tangent to the
orbit
through $\gamma(1)$. We must show that
$J$ is the restriction along $\gamma$
of a Killing field
induced by the $\mathcal P(G,H\times K)$-action.

Let
$\gamma(0)=u\in V$.
Since $\mathcal P(G,H\times K)\supset\Omega_1(G)$,
we
have that the normal space to the orbit through $u$
is contained in the
horizontal subspace $\mathcal H(u)$.
Let $g$ be the unique element
of
$\mathcal P(G,1\times G)$ which satisfies $u=g*\hat 0$.
Now
$\gamma(t)=u+t\Ad{g}\hat X$ for some $X\in\Lg$.
Since $J(0)=0$,
setting
$\gamma_\alpha(t)=u+t\Ad{g}(\hat X+\alpha \hat
Y)=g*\widehat{t(X+\alpha
Y)}$,
where $J'(0)=\Ad{g}\hat Y$ for some
$Y\in\Lg$,
defines a variation of $\gamma=\gamma_0$ through
transversal
geodesics which induces $J$.
Since $\varphi:V\to G$ is an equivariant
Riemannian submersion with
respect to the homomorphism
$h\in\mathcal
P(G,H\times K)\mapsto(h(0),h(1))\in H\times K$,
and each $\gamma_\alpha$ is
a horizontal curve with respect to $\varphi$,
we have that
$\{\varphi\gamma_\alpha\}$ is a variation of
$\bar\gamma=\varphi\gamma$
through transversal geodesics in $G$.
Moreover, the associated
transversal
Jacobi field $\bar J$ along $\bar\gamma$ satisfies $\bar
J(0)=0$,
$\bar J'(0)=d\varphi_u(\Ad{g}\hat Y)=Y\cdot a$,
where
$a=g(1)^{-1}=\varphi(u)$,
and $\bar J(1)$ is tangent to the
orbit
through $\bar\gamma(1)$. It follows by variational completeness of
the
action of $H\times K$ on $G$ that $\bar J$ is the restriction
along
$\bar\gamma$ of an $H\times K$-Killing field on $G$,
namely $\bar J(t)
= Z\cdot \bar\gamma(t)-\bar\gamma(t)\cdot W$
for some $Z\in\Lh$,
$W\in\Lk$.

Note that $\bar\gamma(t)=e^{tX}a$.
Therefore we can write $\bar
J(t)=e^{tX}\cdot (\Ad{e^{-tX}}Z-\Ad{a}W)\cdot
a$. 
It follows that $0=\bar
J(0)=(Z-\Ad a W)\cdot a$, so
that
$W=\Ad{a^{-1}}Z\in\Lk\cap\Ad{a^{-1}}\Lh=\Lq^a$ and
$\bar
J(t)=e^{tX}\cdot (\Ad{e^{-tX}}Z-Z)\cdot a$.
Then we have that $Y\cdot a=\bar
J'(0)=[Z,X]\cdot a$, which gives
$Y=[Z,X]$. 

Finally, consider the
one-parameter subgroup $\{g_\alpha\}$
of $\mathcal P(G,H\times K)$ given
by
$g_\alpha=ge^{\alpha Z}g^{-1}$. Note that $g_\alpha(0)=e^{\alpha Z}\in
H$
and $g_\alpha(1)=a^{-1}e^{\alpha Z}a\in K$,
so that $g_\alpha$ is well
defined. We next show that the Killing field
induced by $\{g_\alpha\}$
coincides with $J$ along $\gamma$.
It suffices to
compare their initial
values.
Since
$g_\alpha*\gamma=(ge^{\alpha Z})*t\hat X$, we have
\[
ge^{\alpha Z}*\hat 0=-(ge^{\alpha Z})'(ge^{\alpha
Z})^{-1}=-g'g^{-1}=g*\hat
0=u, \]
and
\begin{eqnarray*}
 \frac{\partial}{\partial
t}\Big\arrowvert_{t=0}
\frac{\partial}{\partial\alpha}\Big\arrowvert_{\alpha
=0}
ge^{\alpha Z}*t\hat X
&=&
\frac{\partial}{\partial\alpha}\Big\arrowvert_{\alpha=0}
\frac{\partial}
{\partial t}\Big\arrowvert_{t=0}
t\Ad{ge^{\alpha Z}}\hat X + u \\
&=&
\frac{\partial}{\partial\alpha}\Big\arrowvert_{\alpha=0}
\Ad{ge^{\alpha
Z}}\hat X \\
&=& \Ad{g}\widehat{[Z,X]} \\
&=& \Ad{g}\hat
Y.
\end{eqnarray*}
This completes the proof. \EPf

\begin{lem}
If the action
of $\mathcal P(G,H\times K)$ on $V$ is
variationally complete, then it is
hyperpolar.
\end{lem}

\Pf Let $N_0$ in $V$ be an orbit
which is the
preimage of
a principal orbit of $H\times K$ in $G$ and
$\gamma$ a geodesic
starting orthogonally to
$N_0$ in
$p$, i.e.,
$\xi=\gamma'(0)$ is in the
normal space $\nu_p(N_0)$. Let $N_1$ be the orbit
through
$q=\gamma(1)$
and
assume that also
$N_1$ is the preimage of a principal orbit in $G$. It
follows from
variational completeness that
$\gamma(1)$ is not a focal point
of $N_0$ along $\gamma$. We will prove that
the tangent spaces
of $N_0$ at
$\gamma(0)$ and of $N_1$ at $\gamma(1)$ coincide if considered
to be affine
subspaces of
$V$. It follows that the normal spaces of $N_0$ at $\gamma(0)$
and of $N_1$
at $\gamma(1)$ coincide
which obviously implies that the action
is hyperpolar.

Let $E_p\subset T_pN_0$ be the direct sum of the eigenspaces
of the
Weingarten
operator
$A_\xi^{N_0}$ corresponding to the nonvanishing
eigenvalues and define a
corresponding subspace
$E_q\subset T_qN_1$ with
respect to the Weingarten map $A_\xi^{N_1}$.
Similarly let $Z_p$ be the zero
eigenspace in $T_pN_0$ and $Z_q$ the zero
eigenspace
in $T_qN_1$.
We will
first prove that $E_p=E_q$ using an argument from~\cite{DO}. Then we
will
prove that
$Z_p=Z_q$. It will follow that $T_pN_0=T_qN_1$ finishing
the
proof.

Let $X\in E_p$ be an eigenvector corresponding to the
nonvanishing
eigenvalue $\lambda$ of $A_\xi^{N_0}$.
Then $J(t)=(1-\lambda t)X$ is a transversal
Jacobi field along $\gamma$
that is tangent to the orbit $N_0$ and vanishes in $t=\frac{1}{\lambda}$.
Variational completeness now implies that this Jacobi field is induced by
the action.
Therefore
$J(t)$ is tangent to the orbit through $\gamma(t)$ for every $t$.
Since $\gamma(1)$ is not a focal point of $N_0$ along
$\gamma$, we have that $1-\lambda\ne0$ and then
$X$ lies in $T_qN_1$.
It is now clear have that $X$ is
an eigenvalue of $A_\xi^{N_1}$ corresponding to a nonzero eigenvalue.
This proves that $E_p\subset E_q$.
Analogous arguments imply $E_q\subset E_p$ proving
$E_p=E_q$. 
 
It is left to prove that $Z_p=Z_q$.
 Let $X\in Z_p$. Let $\{\psi_\alpha\}$
be a one-parameter subgroup of $\mathcal P(G,H\times K)$ such that
$$X=\frac{d}{d\alpha}\Big|_{\alpha=0}\psi_\alpha p.$$
The variation $\{\psi_\alpha\gamma\}$
induces a Jacobi field $J$ along $\gamma$ with
$J(0)=X$.
We first show that $J(t)=X+t\eta$ where $\eta\in \nu_p(N_0)$. In fact
$$J'(0)= \frac{d}{d\alpha}\Big|_{\alpha=0}{\psi_\alpha}_*\xi=
-A_{\xi}X+\nabla^\perp_{\frac{d}{d\alpha}}
{\psi_\alpha}_*\xi\,|_{\alpha=0}=\nabla^\perp_{\frac{d}{d\alpha}}
{\psi_\alpha}_*\xi\,|_{\alpha=0}$$
which shows that $J'(0)$ lies in $\nu_p(N_0)$ as we wanted to
prove.

Our next goal is to show that $\eta=0$. Assume that $\eta\ne 0$. We know
that $X+\eta\in
T_qN_1$. Then 
$\eta\not\in\nu_q(N_1)$ since $\langle
X+\eta,\eta\rangle=\vert\vert\eta\vert\vert^2\ne 0$. This implies that
$X+\eta\not\in Z_q$ since otherwise $J'(1)=\eta$ would be a normal
vector.  Let
$Y\in E_q$ be the image of
$X+\eta$ under the orthogonal projection along $Z_q$ into $E_q$. We have
that $Y\ne 0$
or, which is the same,
$\langle Y,X+ \eta\rangle\ne 0$.  We know that $\langle Y,X\rangle =0$
since $E_p=E_q$.
We also have that $\langle Y,\eta\rangle=0$ since $Y$ is a tangent vector of
$N_0$ in $p$
and $\eta$ is a normal vector. Hence we also get $\langle Y,X+ \eta\rangle=
0$
which is a
contradiction. This proves that $\eta=0$. We therefore have that $J(t)=X$
which implies
that $X\in T_qN_1$ and it is clear that $X\in Z_q$.
It follows that $Z_p\subset Z_q$. The proof that $Z_q\subset Z_p$ is
analogous. Hence
$Z_p=Z_q$. \EPf

\medskip

With the following lemma we finish the proof of the theorem in
the
introduction.

\begin{lem}
If the action of $\mathcal P(G,H\times K)$ on $V$ is
hyperpolar, then the action of $H$ on $M$ is also
hyperpolar.
\end{lem}

\Pf Let $\Sigma$ be a section of the action of $\mathcal P(G,H\times K)$ on
$V$
and let $A$ be the image of $\Sigma$ under $\varphi$ in $G$. We will show
that 
$A$ is a flat section of the action of $H\times K$ on $G$. Since $\Sigma$ is
horizontal with respect to the Riemannian submersion $\varphi$, we have that
$\varphi|_\Sigma:\Sigma\to G$ is an isometric immersion.
Property (c) before Lemma 2
implies that $d\varphi_u(T_u\Sigma)$ is perpendicular
to the $H\times K$-orbit
through $\varphi(u)$ in $G$ for all $u\in V$. It follows that $A$ is
a submanifold in $G$ which meets all orbits perpendicularly. Moreover $A$ is
flat
since $\varphi|_\Sigma$ is an isometric immersion. This finishes the proof
that 
the action of $H\times K$ on $G$ is hyperpolar. Similar arguments show that
$\pi(A)$ is a flat section of the action of $H$ on $M$. \EPf

\medskip

It is interesting to remark
that it follows from the results in section~2 of~\cite{HPTT} that,
in the case where the metric in $G$ is induced from the Killing form
of $\Lg$, we have that $A$ is a torus, and hence $A$ and $\pi(A)$ are
\emph{properly embedded} sections. In the case where the metric in $G$ is
not induced from the Killing form of $\Lg$, we do not know of
any example of a hyperpolar action with
a section that is not properly embedded. Here we have assumed the
symmetric space $M$ to be of compact type. On the other
hand,  if the symmetric space $M$ is a torus, a hyperpolar action
can have a 
section that is not properly embedded.
To see this let $T^k$ be a proper subtorus of $M$
that we think of as a Lie group. The orbits of $T^k$ are the cosets of
$T^k$,
and the image $\Sigma$ under the exponential map of the normal
space of $T^k$ in $M$ at the identity element
is clearly a section.  If $M$ is not a rational torus,
the subtorus $T^k$ can be chosen so that $\Sigma$ is not
properly embedded, see section~2 in~\cite{HPTT}.

\bibliographystyle{amsplain}
\bibliography{paper}

\providecommand{\bysame}{\leavevmode\hbox to3em{\hrulefill}\thinspace}
\begin{thebibliography}{10}

\bibitem{B}
R.~Bott, \emph{An application of the {M}orse theory to the topology of {L}ie
  groups}, {Bull. Soc. Math. France} \textbf{84} (1956), 251--281.

\bibitem{BS}
R.~Bott and H.~Samelson, \emph{Applications of the theory of {M}orse to
  symmetric spaces}, {Amer. J. Math.} \textbf{80} (1958), 964--1029, correction
  in {Amer. J. Math. {\bf83} (1961), 207--208}.

\bibitem{C}
L.~Conlon, \emph{Variational completeness and {$K$}-transversal domains}, {J.
  Differential Geom.} \textbf{5} (1971), 135--147.

\bibitem{D}
J.~Dadok, \emph{Polar actions induced by actions of compact {L}ie groups},
  {Trans. Amer. Math. Soc.} \textbf{288} (1985), 125--137.

\bibitem{DO}
A.~J. {Di Scala and C. Olmos}, \emph{Variationally complete representations are
  polar}, {Proc.~Amer.~Math.~Soc.} \textbf{129} (2001), 3445--3446.

\bibitem{GTh}
C.~Gorodski and G.~Thorbergsson, \emph{Representations of compact {L}ie groups
  and the osculating spaces of their orbits}, Preprint, University of Cologne,
  2000 (also E-print math.~DG/0203196).

\bibitem{HPTT2}
E.~Heintze, R.~S. Palais, C.-L. Terng, and G.~Thorbergsson, \emph{Hyperpolar
  actions and $k$-flat homogeneous spaces}, {J. Reine Angew. Math.}
  \textbf{454} (1994), 163--179.

\bibitem{HPTT}
\bysame, \emph{Hyperpolar actions on symmetric spaces}, Geometry, Topology, and
  Physics for Raoul Bott (S.~T. Yau, ed.), Conf. Proc. Lecture Notes Geom.
  Topology, {IV}, International Press, Cambridge, MA, 1995, pp.~214--245.

\bibitem{H}
R.~Hermann, \emph{Variational completeness for compact symmetric spaces},
  {Proc. Amer. Math. Soc.} \textbf{11} (1960), 544--546.

\bibitem{K}
A.~Kollross, \emph{A classification of hyperpolar and cohomogeneity one
  actions}, {Trans. Amer. Math. Soc.} \textbf{354} (2002), 571--612.

\bibitem{T1}
C.-L. Terng, \emph{Proper {F}redholm submanifolds of {H}ilbert space}, {J.
  Differential Geom.} \textbf{29} (1989), 9--47.

\bibitem{T2}
\bysame, \emph{Polar actions on {H}ilbert space}, {J. Geom. Anal.} \textbf{5}
  (1995), 129--150.

\bibitem{TTh}
C.-L. Terng and G.~Thorbergsson, \emph{Submanifold geometry in symmetric
  spaces}, { J. Differential Geom.} \textbf{42} (1995), 665--718.

\end{thebibliography}

\bigskip

\parbox[t]{7cm}{\footnotesize\sc Instituto de Matem\'atica e Estat\'\i
stica\\
                Universidade de S\~ao Paulo\\
                Rua do Mat\~ao, 1010\\
                S\~ao Paulo, SP 05508-090\\
                Brazil\\ \hfill\\
                E-mail: {\tt gorodski@ime.usp.br}}{}\hfill{}
\parbox[t]{7cm}{\footnotesize\sc Mathematisches Institut\\
                Universit\"at zu K\"oln\\
                Weyertal 86-90\\
                50931 K\"oln\\
                Germany\\ \hfill\\
                E-mail: {\tt gthorber@mi.uni-koeln.de}}\hfil{}

\end{document}